\documentclass{article}
\usepackage{ifthen}
\usepackage[latin1]{inputenc}
\usepackage{amssymb}
\usepackage[dvips]{graphicx,epsfig}
\usepackage{color} 
\usepackage{lscape,graphics}

\usepackage{ifthen}
\usepackage[latin1]{inputenc}
\usepackage{amsmath,amssymb,amsthm}

\newcommand{\A}{\mathbb{A}}
\def\AA{{\mathbb A}}

\def\CC{{\mathbb C}}

\def\GG{{\mathbb G}}
\def\HH{{\mathbb H}}

\def\NN{{\mathbb N}}

\def\PP{{\mathbb P}}
\def\QQ{{\mathbb Q}}

\def\ZZ{{\mathbb Z}}

\def\fxo{{\ensuremath \mathcal O}}

\newcommand{\fd}{\ensuremath{\rightarrow}}

\newcommand{\findem}{\nolinebreak\vspace{\baselineskip} \hfill\rule{2mm}{2mm}}
\renewcommand{\phi}{\ensuremath{\varphi}}
\newcommand{\inc}{\ensuremath{\subset}}

\newcommand{\ox}{\otimes }
\newcommand{\pla}{\A^2}
\newcommand{\plp}{\PP^2}
\renewcommand{\phi}{\ensuremath{\varphi}}
\newcommand{\x}{\ensuremath{\times}}

\newcommand{\s}{Spec\ }

\newcommand{\codim}{\mbox{codim}\;}

\newtheorem{nt}{Notation}

\newtheorem{prop}[nt]{Proposition}

\newtheorem{exercice}{Exercice} 
\newcounter{numeroquestion} 
\newcounter{numerosousquestion}

\newcommand{\sousquestion}{\ifthenelse{\value{numerosousquestion}=1}{}{\\}\textbf{\roman{numerosousquestion})} \addtocounter{numerosousquestion}{1}}
\newtheorem{correction}{Correction}

\newenvironment{listecompacte}
{\begin{list}
    {\ensuremath{\bullet}}
    {\setlength{\topsep}{2pt}
      \setlength{\itemsep}{1pt} \setlength{\parsep}{0pt}}
}
{\end{list}
}

\newtheorem{coro}[nt]{Corollary} 
\newtheorem{defi}[nt]{Definition} 
\newtheorem{defprop}[nt]{Definition-Proposition} 
 
\newtheorem{ex}[nt]{Example}

\newtheorem{rem}[nt]{Remark} 
\newtheorem{thm}[nt]{Theorem}

\newenvironment{dem}{\noindent\textit{Proof.} }{\findem}

\begin{document}
\sloppy
\title{Intersection theory on punctual Hilbert schemes and graded
  Hilbert schemes}
\author{Laurent Evain (laurent.evain@univ-angers.fr)}
\maketitle


\section*{Abstract}
\label{sec:abstract}

The rational Chow ring 
 $A^*(S^{[n]},\QQ)$ of the Hilbert
scheme $S^{[n]}$ parametrising the length $n$ zero-dimensional subschemes of 
a toric surface $S$ can be described with the help of equivariant techniques. In this paper, 
we explain the general method and we illustrate it through many examples. In the last 
section, we present results on the  
intersection theory of graded Hilbert schemes.

\section*{R\'esum\'e}
Les techniques \'equivariantes permettent de d\'ecrire l'anneau de
Chow rationnel $A^*(S^{[n]},\QQ)$ du sch\'ema de Hilbert $S^{[n]}$
param\'etrant les sous-sch\'emas ponctuels de longueur $n$ d'une
surface torique $S$. Dans cet article, nous pr\'esentons la
d\'emarche g\'en\'erale et nous l'illustrons au travers de nombreux
exemples. La derni\`ere section expose des r\'esultats de
th\'eorie d'intersection sur des sch\'emas de Hilbert gradu\'es. 

\section*{Introduction}
\label{sec:introduction}

Let $S$ be a smooth projective surface and $S^{[n]}$ the Hilbert
scheme parametrising the length $n$ zero-dimensional subschemes of
$S$.  How to describe the cohomology ring 
$H^*(S^{[n]},\QQ)$ and the Chow ring 
$A^*(S^{[n]},\QQ)$ ?

A first approach is based on the work of Nakajima, Grojnowski
and Lehn among others
\cite{lehn_sorger01:cup_product_on_Hilbert_schemes},
\cite{vasserot01:anneauCohomologieHilbert},
\cite{lehn_sorger02:cup_product_on_Hilbert_schemes_for_K3},
\cite{liQinWang04mathAG:cohomoDesSchemaHilbertSurface=FibreSurP1},
\cite{costello-grojnowski03:CohomoSchemaHilbertPonctuel}.
The direct sum $\oplus_{n\in \NN}  H^*(S^{[n]},\QQ)$ 
is an (infinite dimensional) irreducible representation and
carries a Fock space structure \cite{nakajima97:_heisenberg_et_Hilbert_schemes}. 
Lehn settles a connection between the Fock space structure and the 
intersection theory of the Hilbert scheme via the action of the
Chern classes of tautological bundles
\cite{lehn99:_chern_classes_of_tautological_sheaves_on_Hilbert_schemes}.   


An other method, independent of the Fock space formalism introduced by Nakajima, has been 
developed in \cite{evain2007:chowRingHilbertTorique} when $S$ is a toric surface.
The point is that the extra structure coming from the torus action 
brings into the scene an equivariant Chow ring which is easier to 
compute than the classical Chow
ring. The classical Chow ring is a quotient of the 
equivariant Chow ring.

The computations of this equivariant approach are explicit. 
They rely on the standard description 
of the cohomology of the 
Grassmannians and on a description of 
the tangent space to the Hilbert scheme at fixed points.

The main goal of this paper is to present this equivariant
approach. We follow the general theory and we illustrate it with 
the case $S=\PP^2$ and $n=3$ as the main example.

In the last section, we bring our attention to 
graded Hilbert schemes, which played an important role 
in the equivariant computations.  We present results 
on the set theoretic intersection of Schubert cells, which  
suggest that 
intersection theory on graded Hilbert schemes could be
described in terms of combinatorics of plane partitions.

Throughout the paper, we use the formalism of Chow rings and work over any
algebraically closed field $k$.
When $k=\CC$, the Chow ring co{\"\i}ncides with
usual cohomology since the action of the two-dimensional torus $T$ on 
$S$ induces an action of $T$ on $S^{[d]}$ with a finite number of fixed
points.

\section{Equivariant intersection theory}
\label{sec:equiv-inters-theory}

\subsection{General results}
\label{sec:equiv-chow-rings}

In this section, we recall the facts about equivariant Chow rings that we
need. To simplify the presentation, we work with
rational coefficients and the notation $A^*(X):=A^*(X,\QQ)$ denotes 
the rational Chow ring. 

The construction of an equivariant Chow ring associated with an
algebraic space endowed with an action of a linear algebraic group has
been settled by  Edidin and Graham
\cite{edidinGraham:constructionDesChowsEquivariants}.
Their construction is modeled after the Borel construction in equivariant
cohomology.

\begin{prop}
  Let $G$ be an algebraic group,  $X$ an equidimensional 
  quasi-projective scheme with
  a linearized $G$-action and $i,j\in \ZZ$, $i\leq dim(X), j\geq 0$. There exists a
  representation $V$ of $G$ such that
  \begin{listecompacte}\item 
    $V$ contains an open set $U$ on which $G$ acts freely,
  \item $U\fd U/G$ exists as a scheme and is a principal $G$ bundle,
  \item $\codim_V V\setminus U> \dim(X)-i$.
  \end{listecompacte}
  The quotient $X_G=(X\x U)/G$ under the diagonal action exists as a
  scheme. 
  The groups $A_i^G(X):=A_{i+\dim(V)-\dim(G)}(X_G)$ and 
  $A^j_G(X):=A_{\dim(X)-j}^G(X)$ are independent of the choice of the
  couple $(U,V)$.  
\end{prop}

\begin{defi}
  The group $A^i_G(X)$ is by definition the equivariant Chow group of
  $X$ of degree $i$.  
\end{defi}

\begin{ex}
  If $G=T= (k^*)^n$ is a torus, then a possible choice for the couple $(U,V)$ is $V=(k^l)^n$
  with $l>>0$, and $U=(k^l-\{0\})^n$ with  $T$ acting on $V$ by
  $(t_1,\dots,t_n)(x_1,\dots,x_n)=(t_1x_1,\dots,t_nx_n)$. The quotient
  $U/T$ is isomorphic to $(\PP^{l-1})^n$. 
\end{ex}

\begin{ex}
  Let $p$ be a point and $T=(k^*)^n$ the torus acting trivially on
  $p$. Then $A_T^*(p)\simeq \QQ[h_1,\dots,h_n]$ where
  $h_i$ has degree $1$ for all $i$. 
\end{ex}
\begin{dem}
  By the above example, $A^*_T(p)=\lim_{l\fd
    \infty}A^*((\PP^{l-1})^n)=\lim_{l\fd \infty}
    \QQ[h_1,\dots,h_n]/(h_1^{l},\dots,h_n^{l})=\QQ[h_1,\dots,h_n]$, 
    where $h_i$ has degree $1$ ( according to the definition of the
    equivariant Chow group, the limit considered is a degreewise
    stabilisation thus the
    limit is the polynomial ring and not a power series ring). 
\end{dem}

If $X$ is smooth, then $(X\x U)/G$ is smooth too and 
$A_T^*(X)$ is a ring : the intersection of
two classes $u,v\in A_T^*(X)$ takes place in the Chow ring $A^*((X\x U)/G)$. 
\begin{ex}
  The isomorphism $A_T^*(p)\simeq \QQ[h_1,\dots,h_n]$ of the last
  example is an isomorphism of rings. 
\end{ex}

\begin{defi}
  Let $E$ be a $G$-equivariant vector bundle on $X$ and $E_G\fd X_G$ the
  vector bundle with total space $E_G=(E\x U)/G$. The equivariant
  Chern class $c_j^G(E)$ is defined by $c_j^G(E)=c_j(E_G)\in
  A^j(X_G)=A_G^j(X)$. 
\end{defi}

The identification of $A_T^*(p)$ with a ring of polynomials 
$R$ can be made intrinsic using equivariant Chern
classes.
\begin{prop} Let $\hat T$ be the character group of a torus $T\simeq (k^*)^n$. 
  Any character $\chi\in \hat T$ defines a one-dimensional
  representation of $T$ by $t.k=\chi(t)k$, hence an equivariant bundle
  over the point and an equivariant Chern
  class $c^T_1(\chi)$. The map $\chi\fd c_1^T(\chi)\in A^1_T(p)$ extends to an
  isomorphism $R=Sym_\QQ(\hat T) \fd A_T^*(p)$, where  $Sym_\QQ(\hat T)$
 is the symmetric algebra over $\QQ$ of the group $\hat T$. 
\end{prop}

\begin{ex}\label{exemple:projectiveSpace}
  Let $T=k^*$ be the one dimensional torus acting on the projective space
  $\PP^r=Proj \ k[x_0,\dots,x_r]$ by $t.(x_0:\dots:x_r)=(t^{n_0}x_0:\dots:t^{n_r}x_r)$. Then
  $A_T^*(\PP^r)=\QQ[t,h]/p(h,t)$  where $p(h,t)=\sum _{i=0}^r h^{r-i}
e_i(n_0t,\dots, n_rt)$, $e_i$ being the $i$-th elementary symmetric
polynomial.  
\end{ex}
\begin{dem}
  $X_T$ is
  the $\PP^r$ bundle $\PP(\fxo(n_0)\oplus \dots \oplus \fxo(n_r))$
  over $\PP^{l-1}$. The
  rational Chow ring of this projective bundle is $\QQ[h,t]/(p(h,t),t^l)$. We
  have the result when $l$ tends to $\infty$. 
\end{dem}

\begin{ex}\label{exemple:Grassmannian}
  Let $V$ be a representation of $G$ and $G(k,V)$ the corresponding
  Grassmannian. Then $A_G^*(G(k,V))$ is generated as an $R$ module
  by the equivariant Chern classes of the universal quotient bundle. 
\end{ex}
\begin{dem}
  The quotient $(G(k,V)\x U)/G$ is a Grassmann bundle over $U/G$ with
  fiber isomorphic to $G(k,V)$. Since the Chow rings of Grassmann
  bundles are generated over the Chow ring of the base by the Chern classes of the
  universal quotient bundle, the result follows. 
\end{dem}

\subsection{Results specific to the action of tori}
\label{sec:results-spec-acti}

Brion 
\cite{brion97:_equivariant_chow_groups}
pushed further the theory of equivariant Chow rings when the group is a torus 
$T$ acting on a variety $X$. 

\begin{thm}\cite{brion97:_equivariant_chow_groups}\label{thr:restrictionAuxPointsFixesInjective}
  Let $X$ be a smooth projective $T$-variety. 
  The restriction morphism $i_T^*:A_T^*X \fd A_T^*X^T$ is injective. 
\end{thm}

\begin{ex} Let $T=k^*$ act on
$\PP^1$ by $t.(x:y)=(tx:y)$. The inclusion
  $i_T^*:A_T^*(\PP^1)\fd A_T^*(\{0,\infty\})=R^2$ 
  identifies $A_T^*(\PP^1)$ with the couples $(P,Q)$ of polynomials
  $\in R=\QQ[t]$ such that $P(0)=Q(0)$. 
\end{ex}
\begin{dem}
  Let $V=Vect(x,y)$ be the 2 dimensional vector space with
  $\PP(V)=\PP^1$. By the above $A_T^*(\PP^1)$ is generated by the Chern
  classes  of the universal quotient bundle as an $R$-module. On the point
  $\infty = ky\in \PP(V)$, the quotient bundle $Q$ is isomorphic to $kx$ and $T$ acts with
  character $t$. Thus
  $c_1(Q)_{\infty}=t$. Similarly, the restriction of
  $Q$ to the point $0=kx$ is a trivial equivariant bundle
  and $c_1(Q)_{0}=0$. Thus $c_1(Q)$ restricted to $\{0,\infty\}$ is
  $(0,t)$. Obviously, $c_0(Q)=(1,1)$. Thus
  $A_T^*(\PP^1)=\QQ[t](0,t)+\QQ[t](1,1)$ as expected.  
\end{dem}

If $T' \inc T$ is a one codimensional torus, the localisation morphism
$i_T^*$ factorizes: $ A_T^*(X)\fd
A_T^*(X^{T'})\stackrel{i_{T'}^*}{\fd}A_T^*(X^T)=R^{X^T}$. 
Brion has shown 
\begin{thm}\cite{brion97:_equivariant_chow_groups}
\label{thr:sousToresCodim1} 
Let $X$ be a smooth projective variety with an action of
  $T$. 
  The image $Im(i_T^*)$ satisfies $Im(i_T^*)=\cap_{T'}
  Im(i_{T'}^*)$ where the intersection runs over all subtori $T'$ of
  codimension one in $T$.
\end{thm}

An important point is that the equivariant Chow groups determine 
the usual Chow groups. 
The fibers of $X_T\fd U/T$ are isomorphic to $X$. 
Let $j:X \fd X_T$ be the inclusion of a fiber and 
$j^*:A_T^*(X)\fd A^*(X)$ the corresponding restriction. 
\begin{thm} \cite{brion97:_equivariant_chow_groups}
Let $R^+=\hat TR \inc R$ be the set of polynomials with
  positive valuation. 
  The morphism $j^*$ is surjective with kernel $R^+ A_T^*(X)$. 
\end{thm}

\begin{ex}
  $A^*(\PP^1)=(\QQ[t](t,0)+\QQ[t](1,1))/(\QQ[t]^+(t,0)+\QQ[t]^+(1,1))
\simeq \QQ[t]/(t^2)$. The isomorphism sends $(P=\sum p_it^i,Q=\sum q_i
t_i)$ with $p_0=q_0$ to $(p_0,p_1-q_1)$. 
\end{ex}

Finally, we have an equivariant Kunneth formula for the restriction to
fixed points, proved in \cite{evain2007:chowRingHilbertTorique}.
\begin{thm} \label{thm:kunneth}
Let $X$ and $Y$ be smooth projective $T$-varieties with finite set of
  fixed points $X^T$ and $Y^T$. Let $A_T^*(X)\inc R^{X^T}$, $A_T^*(Y)\inc
  R^{Y^T}$, and $A_T^*(X\x Y)\inc R^{X^T\x Y^T}$ the realisation of
  their equivariant Chow rings via localisation to fixed points.  
  The canonical isomorphism $ R^{X^T} \ox_R  R^{Y^T}\simeq R^{X^T\x
    Y^T}$ sends  $A_T^*(X)\ox
  A_T^*(Y)$  to $A_T^*(X\x Y)$.
\end{thm}

\begin{ex}\label{ex:produit}
  Let $T$ be the one dimensional torus acting on $\PP^1\x \PP^1$ by
  $t.([x_1,y_1],[x_2:y_2])= ([tx_1,y_1],[tx_2:y_2])$. 
For each  copy of $\PP^1$,  
$A_T^*(\PP^1)$ is generated  as an $R$-module, 
 by the elements $e=1.\{0\}+1.\{\infty\}=(1,1)$ and
 $f=0.\{0\}+t.\{\infty\}=(0,t)$. 
By the Kunneth formula, $A_T^*(\PP^1\x
 \PP^1)\inc \QQ[t]^4$ is generated by the elements
$(1,1,1,1)=(1,1)\ox
 (1,1)$, $(0,t,0,t)=(1,1)\ox (0,t)$,  $(0,0,t,t)=(0,t)\ox (1,1)$,
$(0,0,0,t^2)=(0,t)\ox (0,t)$ where the coordinates are the
coefficients with respect to the four points $(a,b)\in \PP^1\x \PP^1$ with $a,b\in
\{0,\infty\}$. 
\end{ex}

\subsubsection*{The strategy}
\label{sec:strategy}
Let's sum up the situation. 
The equivariant Chow ring $A_T^*(X)$ satisfies the  usual
functorial properties of a Chow ring: there is an induced pushforward
  for a proper morphism ,  an induced pull-back for a flat morphism, equivariant
  vector bundles have equivariant Chern classes...
When the fixed point set $X^T$ is finite, the computations are
identified with calculations in products of polynomial rings.

Since it is possible to recover the usual Chow ring from the
equivariant Chow ring, the point is to compute the equivariant Chow ring and its
restriction to fixed points. The strategy that will be followed 
in the case of the Hilbert schemes is to use theorem \ref{thr:sousToresCodim1} above. 
It is not obvious a priori 
that the geometry and the equivariant Chow rings of $(S^{[n]})^{T'}$ 
and their restriction to fixed points are 
easier to describe than the equivariant Chow ring of the original variety
$S^{[n]}$. This is precisely the work to be done.

\section{Iarrobino varieties and graded Hilbert scheme}
In our computations of the Chow ring of the Hilbert scheme, 
a central role will be played by the Iarrobino varieties or graded
Hilbert schemes that we introduce now.

Fix a set of dimensions $H=(H_d)_{d \in \NN}$ such that $H_d=0$ for
$d>>0$. 
\begin{defprop}
  The Iarrobino variety $\HH_{hom,H}$ is the set of homogeneous ideals $I=\oplus
  I_d\inc k[x,y]$ of colength $\sum_d H_d$ such that
  $codim(I_d,k[x,y]_d)=H_d$. This is a subvariety of 
$\GG=\prod_{d\ s.t.\ H_d\neq 0} Grass(H_d,k[x,y]_d)$. Moreover, $\HH_{hom,H}$ is empty or
irreducible. 
\end{defprop}
\begin{dem}
  Each vector space $I_d$ corresponds to a
  point in the Grassmannian $Grass(H_d,k[x,y]_d)$ and $I$ corresponds
  to a point in the product $\GG=\prod
  Grass(H_d,k[x,y]_d)$. Accordingly, $\HH_{hom,H}$ is a subvariety  of
  $\GG$.  The irreducibility of 
$\HH_{hom,H}$  is shown in 
\cite{iarrobino77:_punctual_hilbert_schemes_AMS}.
\end{dem}

The Chow ring $A_*(\HH_{hom,H})$ is related to
the Chow ring of $\GG$, as shown by 
King and Walter
\cite{king_walter95:generateurs_anneaux_chow_espace_modules}.

\begin{thm}
  Let $i:\HH_{hom,H}\hookrightarrow \GG$ denote the inclusion.
 The pull back $i^*:A^*(\GG)\fd
  A^*(\HH_{hom,H})$ is surjective.
\end{thm}

There are natural generalisations of the Iarrobino varieties, 
introduced by Haiman and Sturmfels 
\cite{haiman_sturmfels02:multigradedHilbertSchemes}
and called graded Hilbert schemes. 
In our case, the graded Hilbert schemes we are interested in are 
the quasi-homogeneous Hilbert
schemes. 

\begin{defprop}
  Let $weight(x)=a\in \NN$,$weight(y)=b\in \NN$ with $(a,b)\neq
  (0,0)$. Consider the set $\HH_{ab,H}$ of quasi-homogeneous ideals
  $I=\oplus_{d\in \NN} I_d\inc k[x,y]$ with $codim(I_d,k[x,y]_d)=H_d$.
 There is a closed 
embedding $i:\HH_{ab,H}\hookrightarrow \GG=\Pi_{d\in \NN, H_d\neq 0}
Grass(H(d),k[x,y]_d)$. 
\end{defprop}

\begin{rem}
  We could also consider the case $a\in \ZZ$ and/or $b\in \ZZ$. However
  when $ab<0$, $\HH_{ab,H}$ would be empty or a point. Moreover, changing
  $(a,b)$ with $(-a,-b)$ gives an isomorphic graduation.  Consequently,
  any non trivial variety $\HH_{ab,H}$ can be realized with $a\geq 0$
  and $b\geq 0$. We thus consider $a\in \NN$ and $b\in \NN$ without
  loss of generality.
\end{rem}

One wants to extend in this context the results by Iarrobino and
King-Walter. 

Extending Iarrobino's irreducibility result is possible, but not immediate, as
Iarrobino's argument does not extend. 
\begin{thm}\cite{evain04:irreductibiliteDesHilbertGradues}
  The graded Hilbert scheme  $\HH_{ab,H}$ is empty or irreducible. 
\end{thm}
\noindent\textit{Idea of the proof.} 
Since $\HH_{ab,H}$ is smooth as the fixed locus of $(\pla)^{[\sum
  H_d]}$ under 
the action of a one
dimensional torus, irreducibility is equivalent to connectedness. To prove 
connectedness, $\HH_{ab,H}$ admits a stratification where
the cells are the inverse images of the product of Schubert cells on
$\GG$
by the immersion $\HH_{ab,H}\fd \GG$. 
Each cell is an affine space. The cells being  connected spaces,
it suffices to connect together the 
different cells to prove the connectedness of $\HH_{ab,H}$. 
To this aim, one writes down 
explicit flat families over $\PP^1$. These flat families 
correspond to curves drawn on
$\HH_{ab,H}$ that give the link between the different 
cells. 
\findem

The arguments of King and Walter generalise easily to the
quasi-homogeneous case. Moreover, the affine plane $\pla$ is a toric
variety. The action of $T=k^*\x k^*$  on $\pla$ induces an action of
$T$ on $\HH_{ab,H}$. It is possible to generalise the results of
King-Walter to the equivariant setting. 
With minor modifications of their method, one gets the following
theorem. 
\begin{thm}
  The natural restriction morphisms $i^*:A^*(\GG)\fd A^*(\HH_{ab,H})$ 
  and $i_T^*:A_T^*(\GG)\fd A_T^*(\HH_{ab,H})$ are surjective.
\end{thm}

\begin{coro}
  The above theorem induces a description of
  $A_T^*(\HH_{ab,H})$. Indeed, if $j: (\HH_{ab,H})^T\fd \HH_{ab,H}$ is
  the inclusion of the $T$-fixed points, then
  $A^*(\HH_{ab,H})=im(j_T^*)=im(ij)_T^*$
  and the computation of $(ij)_T^*:A_T^*(\GG) \fd A_T^*(\HH_{ab,H}^T)$ follows easily
  from the description of $A_T^*(\GG)$ using equivariant Chern classes.  
\end{coro}

\begin{ex}\label{exemple:calculChowHilbertGradue}
  Let $\HH_{hom,H}$ be the Iarrobino variety parametrising the
  homogeneous ideals in $k[x,y]$ with Hilbert function
  $H=(1,1,0,0,\dots)$. The torus $T=k^*$ acts by $t(x,y)=(tx,y)$. 
  The two fixed points are the ideals $(x^2,y)$
  and $(x,y^2)$. The Iarrobino variety $\HH_{hom,H}$ embeds in 
  the Grassmannian $Grass(1,k[x,y]_1)$
  of one dimensional subspace of linear forms. The universal quotient
  $Q$
  over the Grassmannian is a $T$-bundle. Its restriction 
  over the points $(x^2,y)$ and $(x,y^2)$ is a $T$-representation
  corresponding to the characters  $t \mapsto t$ and 
   $t\mapsto 1$. Thus the first Chern
  class of the universal
  quotient is $(t,0)\in
  R^{\HH_{hom,H}^T}=R^2$. Finally
  $A_T^*\HH_{hom,H}= (c_1^T(Q),c_0^T(Q))=R(t,0)+R(1,1)$.
\end{ex}

\section{Geometry of the fixed locus}
\label{sec:geometry-fixed-locus}

\subsection{Geometry of $(S^{[d]})^T$}
\label{sec:geometry-sdt}

Let $S$ be a smooth projective toric surface. The 2-dimensional torus
$T$ which acts on $S$ acts naturally on $S^{[d]}$. 
According to theorems \ref{thr:restrictionAuxPointsFixesInjective} 
and \ref{thr:sousToresCodim1}, one main
ingredient to describe the equivariant Chow ring $A_T^*(S^{[d]})$ is
to describe the geometry of the fixed loci $(S^{[d]})^T$ and
$(S^{[d]})^{T'}$ of the Hilbert scheme $S^{[d]}$ under the action of the two dimensional
torus $T$ and under the action of any one-dimensional torus $T'\inc T$.

Consider the example $S=\plp=Proj\ k[X,Y,Z]$ and $(\plp)^{[3]}$ the
associated Hilbert scheme. The
action of $T=k^*\x k^*$ on $\PP^2$ 
is $(t_1,t_2).X^aY^bZ^c=(t_1X)^a(t_2Y)^bZ^c$. First, we describe the finite set $((\plp)^{[3]})^T$.  
Obviously, a subscheme $Z\in ((\plp)^{[3]})^T$ has a support 
included in $\{p_1,p_2,p_3\}$ where the $p_i$'s are the toric points
of $\plp$ fixed under $T$. Through each toric point, there are two toric lines 
with local equations $x=0$ and $y=0$. Since $Z$ is $T$-invariant, 
the ideal $I(Z)$ is locally 
generated by monomials $x^\alpha y^\beta$. \includegraphics{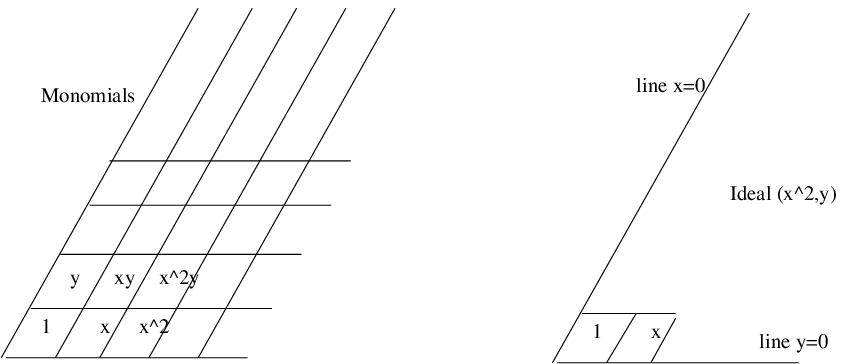}

Using the two lines, we can represent graphically 
the monomials $x^\alpha y^\beta$ as
shown.  An ideal $I\inc k[x,y]$ 
generated by monomials
is represented by the set of monomials which are not in the  
ideal. For instance, the ideal $(x^2,y)$ which does not contain 
the monomials $1,x$ is drawn in the above figure.  

\includegraphics{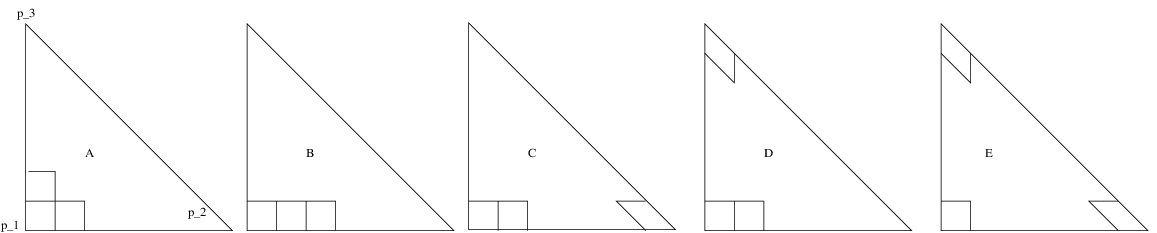}
\begin{prop}
  In $(\plp)^{[3]}$, there are a finite number of subschemes invariant
  under the action of the two-dimensional torus. Up to permutation of
  the projective variables $X,Y,Z$ of $\plp$, there are five such
  subschemes $A,B,C,D,E$. 
\end{prop}
\begin{dem}
  By the above, the invariant subschemes are represented by monomials 
around each toric point of $\plp$. The number of monomials is the
degree of the subscheme, ie. 3 in our situation. Up to permutation of
the axes, all the possible cases $A,B,C,D,E$ are given in the picture.  
\end{dem}

In general, we have:
\begin{prop}
  The points of $(S^{[d]})^{T}$ are in one-to-one correspondence with
  the tuples of staircases $(E_i)_{i\in S^T}$ such that $\sum_{i\in
    S^T}cardinal(E_i)=d$. 
\end{prop}

\subsection{Tangent space at  $p\in (S^{[d]})^{T}$}

We recall the description of the tangent space at a point $p\in  (S^{[d]})^{T}$ 
(\cite{evain04:irreductibiliteDesHilbertGradues}, but see also
\cite{ellingsrud-stromme87:chow_group_of_hilbert_schemes} for 
an other description) .

\includegraphics{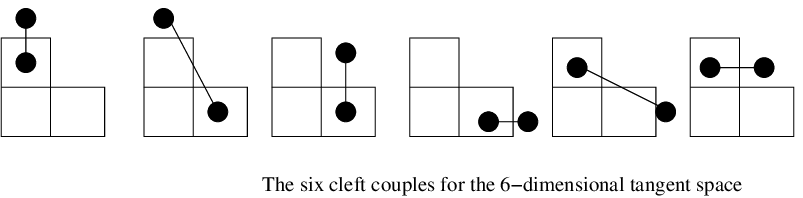}

For simplicity, we suppose that the subscheme $p$ is supported by a single point. 
Recall that we have associated to $p$ the staircase
$F$ of monomials $x^a y^b$ not in $I(p)$ where $x,y$ are the
toric coordinates around the support of $p$.
A cleft for $F$ is a monomial $m=x^ay^b \notin F$ with ($a=0$ or
$x^{a-1}y^b\in F$) and ($b=0$ or
$x^{a}y^{b-1}\in F$). We order the clefts of $F$ according to their
$x$-coordinates: $c_1=y^{b_1},c_2=x^{a_2}y^{b_2},\dots,c_p=x^{a_p}$ 
with $a_1=0<a_2<\dots<a_p$. 
An $x$-cleft couple for $F$ is a couple $C=(c_k,m)$, where $c_k$ is a cleft
($k\neq p$),
$m\in F$, and $mx^{a_{k+1}-a_k}\notin F$. 
By symmetry, there is a notion of
$y$-cleft couple for $F$. The set of cleft
couples is by definition the union of the ($x$ or $y$)-cleft
couples. 
\begin{thm}
  The vector space $T_pS^{[d]}$ is in bijection with the formal sums
  $\sum \lambda_i C_i$, where $C_i$ is a cleft couple for $p$.
\end{thm}

\begin{ex}
  $(\plp)^{[3]}$ is a 6 dimensional variety. A basis for the tangent
  space at a point with local equation $(x^2,xy,y^2)$ is the set of
  cleft couples shown in the figure. 
\end{ex}

With equivariant techniques, it is desirable to describe the tangent
space as a representation. The torus $T$ acts on the
monomials $c_k$ and $m$ with characters $\chi_k$ and $\chi_m$. We
let $\chi_C=\chi_m-\chi_k$. 
\begin{prop}
  Under the correspondence of the above theorem, the cleft couple $C$ is an eigenvector
  for the action of $T$ with character $\chi_C$.
\end{prop}

\subsection{Geometry of $(S^{[d]})^{T'}$}
We come now to the description of $(S^{[d]})^{T'}$ where $T'$ is a one-dimensional subtorus of
$T$.  We start with an example. 

\begin{ex} \label{exempleDeBase}
  If $S=\plp$ and $T'\simeq k^*$ acts
  on $\pla\inc \plp$ via $t.(x,y)=(tx,ty)$ the irreducible components
  of $((\plp)^{[3]})^{T'}$ through the points ${A,B,C,D,E}$ 
  are isomorphic to an isolated point, $ \PP^1,
  \PP^1\x \PP^1,\PP^1\x \PP^1,\PP^2$. 
\end{ex}
\noindent\textit{Proof.}
The tangent space to  $((\plp)^{[3]})^{T'}$ at a point $p\in
((\plp)^{[3]})^{T}$ is the tangent space to  $(\plp)^{[3]}$
invariant under $T'$. Using the description of the tangent space 
to  $(\plp)^{[3]}$ as a representation in the previous section, 
one computes the dimension of the tangent space of $((\plp)^{[3]})^{T'}$
at each of the points  ${A,B,C,D,E}$. The corresponding dimensions are 
$0,1,2,2,2$. 

In particular, an irreducible variety through $A$ (resp. $B,C,D,E$) 
invariant under $T'$ with dimension $0$ (resp. $1,2,2,2$) is the
irreducible component of  $((\plp)^{[3]})^{T'}$ through  ${A}$
(resp. through ${B,C,D,E}$). It thus suffices to exhibit such
irreducible varieties. 

$A$ is an isolated point and there is nothing to do. 

The component $\PP^1$ passing through $B$ can be described
geometrically. The set of lines through the origin of $\AA^2$ is a $\PP^1$. 
To each such line $D$, we consider the subscheme $Z$ of degree $3$
supported by the origin and included in $D$. The set of such subschemes $Z$ 
moves in a $\PP^1$. It is the component of $((\plp)^{[3]})^{T'}$ through $B$. 
This component can be identified with the Iarrobino variety with Hilbert
function $H=(1,1,1,0,0,\dots)$. 

With the same set of lines through the origin, 
one can consider the subschemes $Z$ of
degree $2$ supported by the origin and included in a line $D$. Since $Z$
moves in a $\PP^1$, $Z \cup p$ where $p$ is a point on the line at
infinity moves in a $\PP^1 \x \PP^1$. It is a component through $C$
and $D$. 

Finally, a subscheme  $Z$ of degree $2$ included in the line at
infinity moves in a $\plp$. Thus the union of $Z$ and the origin moves
in a $\plp$ which is the component  of $(S^{[d]})^{T'}$ through $E$. 
\findem

The following proposition says that all but a finite number of representations of
$T'$ give a trivial result. 
\begin{prop}
  Let $a$ and $b$ be coprime with $|a|\geq 3$ or $|b|\geq 3$. Suppose that $T'=k^*$ acts 
  on $\pla\inc \plp$ via $t.(x,y)=(t^ax,t^by)$. Then the irreducible components
  of $((\plp)^{[3]})^{T'}$ through the points ${A,B,C,D,E}$ 
  are isolated points.
\end{prop}
\begin{dem}
  The tangent space at these points is trivial. 
\end{dem}

We see in the examples that  the irreducible components 
of $(\PP^2)^{T'}$ are the three toric points of $\plp$ for a general $T'$.
For some special $T'$, there are two components, namely
a toric point and the line joining the remaining two toric points. 

For a general toric surface $S$, the situation is similar.
\begin{prop}
  For any $T'$ one codimensional subtorus of $T$, $S^{T'}$ is made up
  of isolated toric points $(w_i)$ and of toric lines $(y_j)$ 
joining pairs of the remaining toric points.
\end{prop}

\begin{defi}
  We denote by $PFix(T')=\{w_i\}$ the set of isolated toric points  in
  $S^{T'}$ and
  by $LFix(T')=\{y_j\}$  the set of lines in $S^{T'}$.  
\end{defi}

In the $(\PP^2)^{[3]}$ example, we identified the irreducible components 
of $((\PP^2)^{[3]})^{T'}$ with
products $B_1\x \dots \x B_r$ of projective spaces. 
Some of the projective spaces were identified 
with a graded Hilbert scheme. For instance, 
the component of $((\plp)^{[3]})^{T'}$ through $B$
has been identified with the Iarrobino variety with Hilbert
function $H=(1,1,1,0,0,\dots)$.

In general, the irreducible components 
of $(S^{[d]})^{T'}$ are
products $B_1\x \dots \x B_r$
where the components $B_i$ are projective spaces or graded Hilbert
schemes. 

Let's look at the situation more closely to describe these components.  
If $Z \in (S^{[d]})^{T'}$, the support of $Z$ is invariant
under $T'$. The invariant locus in $S$ is a union 
of isolated points $(w_i)$ and lines $(y_i)$.
We denote by $W_i(Z)$ and $Y_i(Z)$ the subscheme of $Z$ supported 
respectively by the point $w_i$ and by the line $y_i$. By
construction,
we have:
\begin{prop}
A subscheme $Z \in (S^{[d]})^{T'}$ admits a decomposition
$Z=\cup_{w_i \in PFix(T')} W_i(Z) \cup_{y_i \in LFix(T')} Y_i(Z)$. 
\end{prop}
Obviously, $(S^{[d]})^{T'}$ is not irreducible : when $Z$ moves in  
in a connected component, the degree of $W_i(Z)$ and $Y_i(Z)$ should be
constant. But fixing the degree of $W_i(Z)$ and $Y_i(Z)$ is not sufficient 
to characterize the irreducible components of $(S^{[d]})^{T'}$ as shown
by the components identified in example \ref{exempleDeBase}.
\begin{ex} 
  The components of $((\plp)^{[3]})^{T'}$ through $A$ and $B$ are
  2 distinct Iarrobino varieties corresponding to the same degrees $3$
  on the point $(0,0)$ in $\AA^2$ and 0 on the line at infinity.
\end{ex}
The finer 
invariant which distinguishes the irreducible components is 
similar to the one used for the Iarrobino varieties. It is
a Hilbert function taking into account the graduation provided by the
action. 
\begin{ex} Let $O_A=k[x,y]/(x^2,xy,y^2)$ and $O_B=k[x,y]/(y,x^3)$ be
  the ring functions corresponding to the points $A$ and $B$ in
  example \ref{exempleDeBase}. 
  The action of $T'$ on $O_A$ is diagonalizable with
  characters $1,t,t$ whereas the action of $T'$ on $O_B$ acts with
  characters $1,t,t^2$. 
\end{ex}

Let $Z=\s O_Z\in (S^{[d]})^{T'}$. The torus 
$T'$ acts on $O_Z$ with a diagonalizable action.
In symbols, $
O_Z=\oplus V_{\chi}$, where $V_{\chi}\inc O_Z$ is the locus where 
$T'$ acts through the character $\chi \in
      \hat{T'}$.
\begin{defi}
        If $Z \in (S^{[d]})^{T'}$, we define $\underline H_Z:\hat{T'}\fd
        \NN$, $\chi \fd \dim V_{\chi}$ and we let 
        $Z=\cup W_i(Z) \cup Y_i(Z)$ the decomposition of $Z$
        introduced above.
        The tuple of functions $H_Z=(\underline
        H_{W_i(Z)},\underline H_{Y_i(Z)})$
        indexed by the connected components 
        $\{w_i,y_j\}$ of $S^{T'}$ 
        is by definition
        the Hilbert function associated to $Z$.
\end{defi}
\begin{thm}\cite{evain04:irreductibiliteDesHilbertGradues} The set of
  Hilbert functions $H_Z=(\underline H_{W_i(Z)},\underline H_{Y_i(Z)})$ corresponding to the
  subschemes $Z\in (S^{[d]})^{T'}$ is a finite set. Moreover, 
  the irreducible components of $(S^{[d]})^{T'}$ are in one-to-one
  correspondence with this set of Hilbert functions
  $H_Z$.
\end{thm}
The decomposition of  $(S^{[d]})^{T'}$ as a product follows easily from 
the description of the Hilbert functions. Let $B_{w_i}(H)$
be the set of
subschemes $W_i\inc S$, $T'$ fixed, supported by the fixed point $w_i$ with 
$\underline H_{W_i}=H$. Similarly, let $B_{y_i}(H)$ be the set of subschemes 
$Y_i\inc S$, $T'$ fixed, with support in the fixed line $y_i$ and 
$\underline H_{Y_i}=H$. A reformulation of the above is thus:  

\begin{thm}
  $S^{[d],T'}=\cup B_{w_1}(H_{w_1})\x 
\dots B_{w_r}(H_{w_r}) \x B_{y_1}(H_{y_1})\dots \x B_{y_s}(H_{y_s})
$
where $\{w_1,\dots,w_r\}=PFix(T')$ are the isolated fixed points,
$\{y_1,\dots, y_s\}=LFix(T')$ are the fixed lines, and the union runs through 
all the possible tuples of Hilbert functions 
$ (H_{w_1},\dots, H_{w_r},H_{y_1},\dots,H_{y_s})$. 
\end{thm}



The next two propositions identify geometrically the factors
$B_{w_i}(H_{w_i})$ and $B_{y_j}(H_{y_j})$ of the above product. 

By the very definition, we have:

\begin{prop}
For every isolated fixed point $w_i\in PFix(T')$ and every function
$H_{w_i}:\hat T' \fd \NN$,
the variety $B_{w_i}(H_{w_i})$ is a graded Hilbert scheme.   
\end{prop}

\begin{ex}
  In example \ref{exempleDeBase}, take $w=(0,0)\in \pla$ and Hilbert
  function $H(\chi)=1$ for the three characters $\chi=1,t,t^2$ and
  $H(\chi)=0$ otherwise. Then $B_{w}(H)$ is the component of
  $((\plp)^{[3]})^{T'}$ passing through $B$, which has been identified
  with the homogeneous Hilbert scheme $H_{hom,H}$. 
\end{ex}

As for the other components, we have:
\begin{prop}
For every fixed line $y_i\in LFix(T')$ and every function
$H_{y_i}:\hat T' \fd \NN$,
$B_{y_i}(H_{y_i})$ is a product of projective spaces. 
\end{prop}
\noindent
\textit{Illustration of the last proposition on an example. }
Consider $T'=k^*$ acting on  
$\pla\inc \plp$ via $t.(x,y)=(tx,y)$.
The line $x=0$ is $T'$-fixed. We say that a scheme $Z$ is horizontal
of length $n$ if it is in the affine plane and $I(Z)=(y-a,x^n)$, or if it is 
a limit of such schemes when the support $(0,a)$ moves to infinity. 
Two horizontal schemes $Z_1$ and $Z_2$ 
of respective length $n_1\neq n_2$ move in a $\PP^1\x
\PP^1$. When the supports of $Z_1$ and $Z_2$ are distinct, 
$Z_1\cup Z_2$ is a scheme of length $n_1+n_2$. We
thus obtain a rational function $\PP^1\x \PP^1\ \cdots> (\plp)^{[n_1+n_2]}$. 
The schemes being horizontal, the limit of $Z_1 \cup Z_2$ is 
completely determined by the support when the
schemes $Z_1$ and $Z_2$ collide. 
More formally, the rational function extends to a well defined
morphism $\phi:\PP^1\x \PP^1 \fd (\plp)^{[n_1+n_2]}$. This morphism is an
embedding and gives an isomorphism between $\PP^1\x \PP^1$ and 
one of the components $B_{y_i}(H_{y_i})$ introduced above. 

If $n_1=n_2$ in the above paragraph, $\phi$ is not an embedding any
more because of the action of the symmetric group which exchanges 
the roles of $Z_1$ and $Z_2$. But taking the quotient, we obtain an 
embedding $\PP^2=(\PP^1\x \PP^1)/\sigma_2 \fd (\PP^2)^{[n_1+n_2]}$
which is an isomorphism on a component. \findem

\section{Description of the Chow ring}
\label{sec:conclusion}

\subsection{Description using generators}
\label{sec:descr-using-gener}



Let $T'\inc T$ be a one dimensional subtorus. 
Recall that each irreducible component $C$ of $(S^{[d]})^{T'}$
is 
a  product $B_{w_1}(H_{w_1})\x 
\dots B_{w_r}(H_{w_r}) \x B_{y_1}(H_{y_1})\dots \x B_{y_s}(H_{y_s})
$ where $\{w_1,\dots,w_r\}=PFix(T')$ and 
$\{y_1,\dots, y_s\}=LFix(T')$.

The factors $B_{w_i}(H_{w_i})$ are graded Hilbert schemes. We have
seen in example \ref{exemple:calculChowHilbertGradue} the
computation of generators for their equivariant Chow ring. We denote by
$M_{w_i,T',H_{w_i}}$ this equivariant Chow ring. 

The factors $B_{y_i}(H_{y_i})$ are product of projective spaces. We have seen in
examples \ref{exemple:projectiveSpace},
\ref{exemple:Grassmannian} or \ref{ex:produit} the computation of generators for their
equivariant Chow ring.  We denote by
$N_{y_i,T',H_{y_i}}$ this equivariant Chow ring. 

Then the equivariant Chow ring of the component $C$ is given by the
Kunneth formula of theorem \ref{thm:kunneth}: 
\begin{displaymath}
A_T^*(C)=\bigotimes_{w_i\in PFix(T')} M_{w_i,T',H_{w_i}}
\bigotimes_{y_i\in Lfix(T')}N_{y_i,T',H_{y_i}}
\end{displaymath}

When 
$H=(H_{w_1},\dots,H_{w_r},H_{y_1},\dots,H_{y_s})$
runs through the possible Hilbert functions
to describe all the irreducible
components $C$ of $(S^{[d]})^{T'}$ and using theorem \ref{thr:sousToresCodim1},
we finally get: 
\begin{thm}\cite{evain2007:chowRingHilbertTorique}
$$  A_T^*(S^{[d]})=\bigcap_{T'\subset T}\ \bigoplus_{
  \# H=d}
\\
(\ \ \bigotimes_{w_i\in
  PFix(T')}M_{w_i,T',H_{w_i}}  \bigotimes_{y_i\in
  LFix(T')}N_{y_i,T',H_{y_i} \ )}$$
\end{thm}

\subsection{Second description of the Chow ring: From generators to relations}
\label{sec:from-gener-relat}

In the last formula, the modules $M_{w_i,T',H_{w_i}}$ and
$N_{y_i,T',H_{y_i} }$ were described with explicit generators.  It is
possible to adopt the relations point of view rather than the
generators point of view. This is an algebraic trick which relies 
on Bott's integration formula, proved by Edidin and Graham in an
algebraic context. The equivariant Chow ring is then described 
as a set
of tuples of polynomials satisfying congruence relations. 

The proposition below that makes the transition 
from generators to relations
involves equivariant Chern classes of the
restrictions $T_{S^{[d]},p}$ of the tangent bundle $T_{S^{[d]}}$ at
fixed points $p\in (S^{[d]})^T$. Since we have described the fiber of the tangent bundle 
at these points as a $T$-representation, computing the equivariant
Chern classes is straightforward and the set of relations can be
computed.  

\begin{prop}\cite{evain2007:chowRingHilbertTorique}
  Let $\beta_i=(\beta_{ip})_{p\in (S^{[d]})^T}$ be a set of generators of the $\QQ[t_1,t_2]$-module
  $i_T^*A_T^*(S^{[d]})\subset \QQ[t_1,t_2]^{(S^{[d]})^T}$. Let 
$\alpha=(\alpha_p)\in \QQ[t_1,t_2]^{(S^{[d]})^T}$. Then the
  following conditions are equivalent.
  \begin{listecompacte}
  \item 
    $\alpha \in i_T^*A_T^*(S^{[d]}) $
  \item $\forall i$, the congruence 
\begin{displaymath}
\ \sum_{p \in (S^{[d]})^T} (\alpha_p\beta_{ip} \prod_{q\neq
  p}c_{\dim S^{[d]}}^T(T_{S^{[d]},q}))\equiv 0\ (mod\ \prod_{p \in
  (S^{[d]})^T}c_{\dim S^{[d]}}^T(T_{S^{[d]},p}))
\end{displaymath}
holds. 
\end{listecompacte}
\end{prop}

We apply the method to $(\PP^2)^{[3]}$. There are 22 fixed points thus the
equivariant Chow ring is a subring of $\QQ[t_1,t_2]^{22}$. Five of the
fixed points $A,\dots,E$ have been introduced in the
examples. The other fixed points are obtained from these five by a
symmetry. For instance, $A_{12}=\sigma(A)$ where $\sigma$ is the
toric automorphism of $\plp$ exchanging the points $p_1$ and $p_2$.

\begin{thm}\label{thr:leCasHilbTroisP2} \cite{evain2007:chowRingHilbertTorique}
  The equivariant Chow ring  $A_T^*((\PP^{2})^{[3]})\inc
  \QQ[t_1,t_2]^{\{A,A_{12},\dots,E\}}$ 
is the set of linear
  combinations $aA+a_{12}A_{12}+\dots +eE$ satisfying the relations
  \begin{listecompacte}
\item 
$a+a_{13}-d-d_{13}\equiv 0\ (mod \ t_2^2)$
\item $d-d_{13}\equiv 0\ (mod \ t_2)$
\item $a-a_{13}\equiv 0\ (mod \ t_2)$
\item
$a-b\equiv 0\ (mod \ 2t_1-t_2)$
\item
$b-b_{13}\equiv 0\ (mod\ t_2)$
\item
$-b+3c-3c_{12}+b_{12}\equiv 0\ (mod \ t_1^3)$
\item $-b+c+c_{12}-b_{12}\equiv 0\ (mod \ t_1^2)$
\item $3b-c+c_{12}+-3b_{12}\equiv 0\ (mod \ t_1)$
\item
$b-b_{23}\equiv 0\ (mod \ t_2-t_1)$
\item
 $ c-d+c_{23}-d_{23} \equiv 0\ (mod \ (t_1-t_2)^2)$
\item   $c+d-c_{23}-d_{23} \equiv 0\ ( mod\  (t_1-t_2))$
\item $c_{23}-d_{23} \equiv 0\ (mod \ (t_1-t_2))$
\item
$c-c_{13}\equiv 0\ (mod \ t_2)$
\item
$d-2e+d_{12}\equiv 0\ (mod \ t_1^2)$
\item
$d-d_{12}\equiv 0\ (mod \ t_1)$
\item all relations deduced from the above by the action of the
  symmetric group $S_3$. 
  \end{listecompacte}
\item 
The Chow ring $A^*((\PP^{2})^{[3]})$ is the quotient of $A_T^*((\PP^{2})^{[3]})$
by the ideal generated by the elements $fA+\dots +fE$, $f\in \QQ[t_1,t_2]^+$. 
\end{thm}



\section{Graded Hilbert schemes revisited}
\label{sec:homog-vari}

The quasi-homogeneous Hilbert schemes played a central role in the 
computation of $A_T^*(S^{[d]})$ and their Chow ring was computed using equivariant
techniques. In this section, we present a result suggesting
that their Chow ring could admit an alternative description in terms of
combinatorics of partitions. 

To simplify the notations, we restrict from now on our attention 
to the homogeneous case of 
Iarrobino varieties, but the statements below can be formulated in the
quasi-homogeneous case \cite{evain00:cnincidence_cellules_schubert}.


Recall that the graded Hilbert scheme $\HH_{hom,H}$ embeds in a product 
of Grassmannians:  $\HH_{hom,H}\hookrightarrow 
\GG=\prod_{d\in \NN, H_d\neq 0} Grass(H_d,k[x,y]_d)$.
The Grassmannians are stratified by their Schubert cells, constructed
with respect to the flag $F_1=Vect(x^d)\inc F_2=Vect(x^d,x^{d-1}y)
\inc \dots F_{d+1}=k[x,y]_d$. 
The product $\GG$ is stratified by the product of Schubert cells, 
and $\HH_{hom,H}$ is stratified by the restrictions  of these
products of Schubert cells. We still call these 
restrictions Schubert cells on $\HH_{hom,H}$. 
\begin{ex}
  A Schubert cell of a Iarrobino variety $\HH_{hom,H}$ 
contains a unique monomial ideal $I\inc k[x,y]$ that we represent as usual 
by the set of monomials $E(I)=\{x^ay^b \notin I\}$. 
\end{ex}

The Grassmannians in the product $\GG$ are trivial when
$H_d=\dim k[x,y]_d$. When the numbers of non trivial
Grassmannians in $\GG$ is one,  
the inclusion $\HH_{hom,H}\inc \GG$ is an isomorphism.  In this
Grassmannian case, the closures of the Schubert cells
form a base of $A^*(\HH_{hom,H})$ and the intersection 
is classically described in terms of combinatorics involving the partitions associated to
the cells. 
\\[\baselineskip]
\textbf{Question}: \textit{Is it possible to describe the intersection theory in terms of
  partitions when $\HH_{hom,H}$ is not a Grassmannian ? 
}
\\[\baselineskip]

The intersection theory when $\HH_{hom,H}$ is not a Grassmannian is 
more complicated than in the
Grassmannian case. On the set theoretical level, the intersection $\overline C \cap
\overline D$ between the closures of two cells $C$ and $D$ is
difficult to determine. 
The closure $\overline C$ of
a cell $C$ is not a union of cells any more.

However, a
necessary condition for the incidence $\overline C\cap D \neq
\emptyset$ is known and expressed in terms of reverse plane partitions
with shape $E(J)$, where $J$ is the unique monomial ideal in $D$.  

\center{\includegraphics{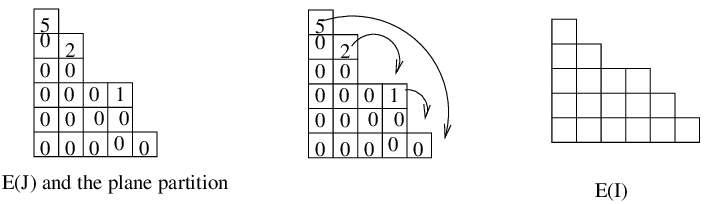}}

\begin{defi}
  A reverse plane partition with shape $E\inc \NN^2$ 
  is a two-dimensional array of integers $n_{ij}, \ (i,j)\in E$
   such that $n_{i,j}\leq n_{i,j+1}$, 	 
$n_{i,j}\leq n_{i+1,j}$ 	
\end{defi}

\begin{defi}
  A monomial ideal $I$ is linked to a monomial ideal $J$  by a reverse plane
  partition $n_{ij}$ with support $E(J)$ if
  $E(I)=\{x^{a+n_{a,b}}y^{b-n_{ab}}, (a,b)\in E(J)\}$
\end{defi}

\begin{ex}
  In the above figure, $E(I)$ is linked to $E(J)$. 
\end{ex}

\begin{defi}
  If $I\inc k[x,y]$ is a monomial ideal of colength $n$, the complement of $I$ is
  the ideal $C(I)$ such that $E(C(I))$ contains the monomials $x^ay^b$, $a< n$,
  $b< n$ with $x^{n-1-a}y^{n-1-b}\in I$ (see the figure below). 
\end{defi}
\includegraphics{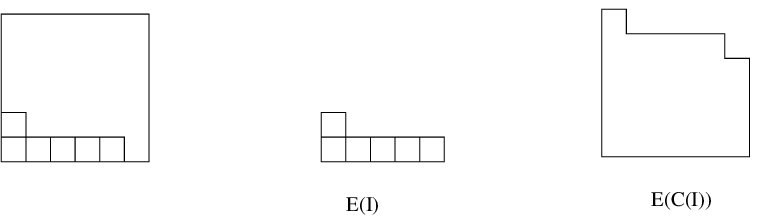}

\begin{thm}\cite{evain00:cnincidence_cellules_schubert}
  Let $C$ and $C'\inc \HH_{hom,H}$ be two cells containing the monomial ideals $I$ and
  $I'$. If $\overline C\cap C'\neq \emptyset$, then
  \begin{listecompacte}
  \item $I$ is
  linked to $I'$ by a reverse plane partition.
  \item $C(I)$ is
  linked to $C(I')$ by a reverse plane partition. 
  \end{listecompacte}
\end{thm}

By analogy with the Grassmannian case, we are led to the following
open question: Can we describe the intersection theory on the Graded
Hilbert schemes in terms of combinatorics of the plane partitions ?

\bibliographystyle{plain} 

\end{document}